\documentclass[10pt,reqno]{amsart}
\usepackage{amsmath,amssymb,mathtools}
\usepackage{microtype}
\usepackage[hidelinks,bookmarksnumbered]{hyperref}
\hypersetup{pdftitle={Algebraic derivations and the first L2-cohomology of free group factors},pdfsubject={Algebraic derivations on free group factors},pdfkeywords={free group factors, algebraic derivations, L2-cohomology}}
\numberwithin{equation}{section}
\newtheorem{theorem}{Theorem}[section]
\newtheorem{proposition}[theorem]{Proposition}
\newtheorem{lemma}[theorem]{Lemma}
\newtheorem{corollary}[theorem]{Corollary}
\theoremstyle{remark}

\newtheorem{example}[theorem]{Example}
\DeclareMathOperator{\Der}{Der}
\DeclareMathOperator{\Inn}{Inn}
\DeclareMathOperator{\Aff}{Aff}
\DeclareMathOperator{\Pack}{Pack}
\DeclareMathOperator{\rk}{rank}
\DeclareMathOperator{\Tr}{Tr}
\DeclareMathOperator{\tr}{tr}
\newcommand{\C}{\mathbb C}
\newcommand{\R}{\mathbb R}
\newcommand{\F}{\mathbb F}
\newcommand{\op}{\mathrm{op}}
\newcommand{\sa}{\mathrm{sa}}
\newcommand{\alg}{\mathrm{alg}}
\newcommand{\ev}{\operatorname{ev}}
\newcommand{\Id}{\mathbf 1}

\author[J. Huang]{Jinghao Huang}
\address{Institute for Advanced Study in Mathematics, Harbin Institute of Technology,
Harbin, 150001, China}
\email{jinghao.huang@hit.edu.cn}

\author[K. Kudaybergenov]{Karimbergen Kudaybergenov}
\address{Institute for Advanced Study in Mathematics, Harbin Institute of Technology,
Harbin, 150001, China and
Suzhou Research Institute of Harbin Institute of Technology, Suzhou, 215104, China}
\email{kudaybergenovkk@gmail.com}

\author[F. Sukochev]{Fedor Sukochev}
\address{School of Mathematics and Statistics, University of New South Wales,
Kensington, 2052, Australia}
\email{f.sukochev@unsw.edu.au}

\author[B. Yan]{Bing Yan}
\address{Institute for Advanced Study in Mathematics, Harbin Institute of Technologies,
Harbin, 150001, China}
\email{bingyan202510@163.com}

\title[First \(L^2\)-cohomology of free group factors]{Algebraic derivations and the first \(L^2\)-cohomology of free group factors}
\date{November 24, 2026}
\subjclass[2020]{Primary 46L10; Secondary 46L54, 46M18}
\keywords{Free group factor, affiliated operator, algebraic derivation, \(L^2\)-Betti number, random matrix}

\begin{document}

\begin{abstract}
Let \(M\) be the von Neumann algebra of a nonabelian free group of finite rank.
We prove that its first \(L^2\)-Betti number in the sense of Connes and
Shlyakhtenko is infinite. We show that every assignment of values in
\(\Aff(M\bar\otimes M^{\op})\) to the standard semicircular generators
extends to an algebraic derivation on all of \(M\). The extension follows
from a rank identity for the Jacobians of finite systems of polynomial
relations. We also prove that, for any \(\mathrm{II}_1\) factor with
separable predual, the existence of an outer derivation with these
coefficients forces the first \(L^2\)-Betti number to be infinite.
\end{abstract}

\maketitle

\section{Introduction}\label{sec:intro}

The study of \(L^2\)-Betti numbers goes back to Atiyah's work on elliptic
operators on normal coverings \cite{Atiyah}. In 2005, Connes and
Shlyakhtenko defined these invariants for tracial \(*\)-algebras using
Hochschild homology \cite{CS}. For the group algebra of a discrete group
with its canonical trace, their definition recovers the \(L^2\)-Betti
numbers of the group \cite[Proposition~2.3]{CS}. In particular, the first
\(L^2\)-Betti number of \(\C\F_n\) is \(n-1\). For its weak closure
\(L(\F_n)\), Connes and Shlyakhtenko's estimates involving free entropy
dimension suggested that \(n-1\) should be a lower bound
\cite[Introduction]{CS}. We prove that the first \(L^2\)-Betti number
of \(L(\F_n)\) is infinite for every integer \(n\ge2\).

Our approach uses Thom's description of the first \(L^2\)-cohomology
in terms of derivations \cite{Thom}. Let \((M,\tau)\) be a finite von
Neumann algebra with a faithful normal tracial state. Put
\[
 N=M\bar\otimes M^{\op},\qquad \tau_N=\tau\otimes\tau^{\op},\qquad
 E=\Aff(N),
\]
where \(M^{\op}\) is the opposite algebra and \(E\) is the algebra of
operators affiliated with \(N\). For \(a\in M\), write
\(L_a=a\otimes1\) and \(R_a=1\otimes a^{\op}\). The bimodule action
on \(E\) is given by \(a\cdot\xi\cdot b=L_aR_b\xi\).
Denote by \(\Der_{\C}(M,E)\) the right \(E\)-module of complex-linear
maps \(\delta:M\to E\) satisfying the Leibniz identity
\begin{equation}\label{eq:leibniz}
 \delta(ab)=L_a\delta(b)+R_b\delta(a),\qquad a,b\in M.
\end{equation}
Its submodule \(\Inn(M,E)\) consists of the inner derivations
\(\delta_\xi(a)=(L_a-R_a)\xi\), with \(\xi\in E\).
In this notation, Thom's duality \cite[Corollaries~3.4 and~3.7]{Thom}
gives
\begin{equation}\label{eq:duality}
 H^1_{\alg}(M,E)=\Der_{\C}(M,E)/\Inn(M,E),\qquad
 \beta_1^{(2)}(M,\tau)=\dim_N H^1_{\alg}(M,E),
\end{equation}
where \(\dim_N\) is the extended von Neumann dimension.

Thom also considered the cohomology obtained by restricting to derivations
that are continuous from the operator norm on \(M\) to the measure
topology on \(E\) \cite{Thom}. Alekseev proved that this continuous
cohomology vanishes for free group factors \cite{Alekseev}, and Popa and
Vaes established the same conclusion for every \(\mathrm{II}_1\) factor
\cite[Theorem~1]{PV}. In particular, every outer derivation constructed
below is discontinuous for these topologies.

Let \(S=(s_1,\ldots,s_n)\) be a standard free semicircular family and
let \(M=W^*(S)\). Each \(W^*(s_j)\) has a Haar unitary generator, and
these unitaries are free, so \(M\) is isomorphic to \(L(\F_n)\).
Since the \(s_j\) are algebraically free, prescribing their images in
\(E\) determines a derivation on the polynomial algebra
\(\C\langle s_1,\ldots,s_n\rangle\). The question is whether this
derivation extends to the weak closure.

\begin{theorem}\label{thm:main}
Let \(2\le n<\infty\), let \(S=(s_1,\ldots,s_n)\) be a standard free
semicircular family, and let \(M=W^*(S)\). For every
\(\eta_1,\ldots,\eta_n\in E\), there exists a derivation
\(\delta\in\Der_{\C}(M,E)\) such that \(\delta(s_j)=\eta_j\) for
\(1\le j\le n\).
\end{theorem}

In other words, evaluation at \(S\) gives a right
\(E\)-linear map
\(\ev_S(\delta)=(\delta(s_1),\ldots,\delta(s_n))^{\mathsf T}\)
from \(\Der_{\C}(M,E)\) onto \(E^n\).
In particular, choose \(\delta(s_1)=0\) and \(\delta(s_2)=1_E\).
Since the spectral distribution of \(s_1\) is diffuse,
\(L_{s_1}-R_{s_1}\) is invertible in \(E\). An inner derivation
vanishing at \(s_1\) therefore vanishes identically, so this choice gives
an outer derivation. To deduce the value of the Betti number, we prove
in Theorem~\ref{dim:dichotomy} that for a \(\mathrm{II}_1\) factor with
separable predual, finiteness of the first \(L^2\)-Betti number forces
every algebraic derivation with coefficients in \(E\) to be inner.

\begin{corollary}\label{cor:main}
Let \(M=L(\F_n)\), where \(2\le n<\infty\). Then
\(H^1_{\alg}(M,E)\ne0\) and \(\beta_1^{(2)}(M,\tau)=\infty\).
\end{corollary}

The main step in the extension theorem concerns relations involving
elements of the weak closure. Let \(Y\) be a finite self-adjoint tuple
in \(M\), and let \(P(S,Y)=0\) be a finite system of self-adjoint
polynomial relations. Write its evaluated Jacobian as \([\,A\ B\,]\),
where the columns of \(A\) correspond to \(S\) and those of \(B\) to
\(Y\). Connes and Shlyakhtenko obtained a relative dimension formula
under a Hilbert--Schmidt commutator condition on the auxiliary elements
\cite[Theorem~4.11]{CS}. For an arbitrary tuple \(Y\) as above, we prove
in Proposition~\ref{prop:rank} that
\[
 \rk_N[\,A\ B\,]=\rk_N B.
\]
Thus every prescribed column \(\eta\in E^n\) admits a solution \(\zeta\)
of the differentiated system \(A\eta+B\zeta=0\).
The rank identity follows by applying the packing estimate of
Section~\ref{sec:geometry} for polynomial gradients to random matrix
approximations of \((S,Y)\). Self-injectivity of \(E\) then gives a
single derivation on \(M\) with the prescribed values.

The same packing estimate also applies to a finite self-adjoint generating
tuple \(X\) admitting matricial microstates. In Section~\ref{sec:entropy}
we show, in particular, that if \(W^*(X)\) is a \(\mathrm{II}_1\) factor
and the modified free entropy dimension \(\delta_0(X)\) exceeds one,
then its first \(L^2\)-Betti number is infinite. Appendix~\ref{q:section}
gives a sharp logarithmic integrability estimate for the solution
matrix \(B^\dagger A\), where \(B^\dagger\) is the Moore--Penrose inverse
of \(B\). Appendix~\ref{cx:sec} gives a complex-geometric proof of the
packing estimate when the auxiliary gradient is of order \(t^2\) at
packing scale \(t\).

\section{A packing estimate for polynomial gradients}\label{sec:geometry}
For a subset \(H\) of a Euclidean space and a number \(t>0\), let \(\Pack_t(H)\) be the supremum of the cardinalities of finite subsets of \(H\) whose distinct points have distance at least \(t\). Let \(\omega_j\) denote the volume of the Euclidean unit ball in \(\R^j\), with \(\omega_0=1\). All finite direct sums of Euclidean spaces carry the square-sum norm.

\begin{lemma}\label{lem:linearpack}
Let \(a\ge1\), \(m\ge0\), and \(d=a+m\). Let \(\Phi:\R^a\times\R^m\to\R\) be a polynomial of degree at most \(p\). Given \(R\ge1\) and \(0<t<1\), let \(\mathcal P_t\) consist of the pairs \((x,z)\in\R^a\times\R^a\) for which there is a vector \(u\in\R^m\) satisfying
\[
 z=\nabla_x\Phi(x,u),\qquad
 |x|,|z|,|u|\le R,\qquad |\nabla_u\Phi(x,u)|\le Rt .
\]
There exists a constant \(C=C(p,R)\ge1\) such that
\begin{equation}\label{eq:packing}
 \Pack_t(\mathcal P_t)\le C^{\,d}t^{-a}.
\end{equation}
\end{lemma}

\begin{proof}
Write \(v=\nabla_u\Phi(x,u)\). In the space of position and momentum coordinates, let \(H\) denote the \(2a\) coordinates associated with \(x,z\), and \(K\) the \(2m\) coordinates associated with \(u,v\). Write \(\Omega_H\) and \(\Omega_K\) for the standard symplectic forms on the two blocks. The symmetry of the Hessian of \(\Phi\) gives
\[
 \sum_{i=1}^a dx_i\wedge dz_i+\sum_{j=1}^m du_j\wedge dv_j=0
\]
on the gradient graph. Rescale the graph by sending \((x,z,u,v)\) to \((x/t,z/t,u,v/t)\). In position--momentum order the rescaled graph \(G\) has the parametrization \(P(q)=(q,g(q))\), where \(q=(x/t,u)\) and \(g\) is a polynomial map of degree at most \(h=\max(1,p-1)\). The pullbacks of the two symplectic forms to \(G\) satisfy
\begin{equation}\label{eq:forms}
 P^*\Omega_H=-t^{-1}P^*\Omega_K .
\end{equation}

Fix a finite \(t\)-separated subset of \(\mathcal P_t\) with \(b\) points. Lift each point to \(G\) using one of its witnesses \(u\). The \(H\)-coordinates of the lifted points have mutual distances at least \(1\). Thus the \(2d\)-dimensional balls of radius \(r=1/3\) about the lifted points are disjoint. If \(v_0\) belongs to one of these balls, a nearest point \(P(q)\) to \(v_0\) exists on the closed graph \(G\). The displacement \(W=v_0-P(q)\) is orthogonal to the tangent space of \(G\) at \(P(q)\). Hence \(W=W(q,c)=(-Dg(q)^{\mathsf T}c,c)\) for some \(c\in\R^d\). Let \(F(q,c)=P(q)+W(q,c)\), and put \(D_0=3(R+1)\). Each nearest point \(P(q)\) is within \(2r\) of the center of the corresponding ball. It follows that the union of the \(b\) balls lies in \(F(U)\), where
\[
 U=\{(q,c): |P_H(q)|\le D_0/t,\ |P_K(q)|\le D_0,\ |W(q,c)|\le r\}.
\]
The set \(U\) is compact since \(q\) and \(c\) are subvectors of \(P(q)\) and \(W(q,c)\), respectively.

The volume of \(F(U)\) is bounded by the integral over \(U\) of the absolute value of the Jacobian of \(F\), by the area formula \cite{Federer}. Up to the choice of orientation, its Jacobian form is
\begin{equation}\label{eq:jacform}
 \frac{(F^*\Omega_H)^a}{a!}\wedge\frac{(F^*\Omega_K)^m}{m!}.
\end{equation}
Consider first a pair of \(H\)-coordinates in \(F=P+W\). Apart from the term involving two \(P\)-differentials, the other three terms in the expansion of the area form contain at most one \(P_H\)-differential each. Replace the sum of the terms involving two \(P_H\)-differentials by \(-t^{-1}P^*\Omega_K\), using \eqref{eq:forms}. The first factor in \eqref{eq:jacform} is now a divided power of a sum of \(m+3a\) simple two-forms. The second factor is a divided power of a sum of \(4m\) simple two-forms involving only \(P_K\) and \(W\). Simple two-forms square to zero and two-forms commute. Consequently the expanded product has at most \(2^{m+3a}2^{4m}\) terms. A term using \(l\) replacements has a scalar factor \(t^{-l}\) and involves \(j\le a-l\) differentials of \(P_H\).

After removing the sign and the factor \(t^{-l}\), each surviving term is the Jacobian form of a map \(\pi:\R^{2d}\to\R^{2d}\) obtained by selecting \(2d\) different coordinates from \((P,W)\). The map \(\pi\) is polynomial of degree at most \(h\). At a regular preimage the complexified fiber has an isolated point, so the affine B\'ezout bound gives at most \(h^{2d}\) regular preimages of any value \cite[Example~8.4.6]{Bezout}. The image \(\pi(U)\) lies in a product of balls of dimensions \(j\) and \(2d-j\), with radii \(D_0/t\) and \(2D_0\), respectively. Since \(\omega_j\omega_{2d-j}\le2^d\omega_{2d}\), a second application of the area formula gives
\[
 t^{-l}\int_U|\det D\pi|
 \le h^{2d}2^d\omega_{2d}(2D_0)^{2d}t^{-l-j}
 \le C_1^{\,d}\omega_{2d}t^{-a}.
\]
Here \(C_1\) depends only on \(p\) and \(R\). Summing over the terms in \eqref{eq:jacform} and dividing by the single-ball volume \(\omega_{2d}r^{2d}\) proves \eqref{eq:packing}.
\end{proof}

\begin{example}
The scale \(t\) in the bound on \(\nabla_u\Phi\) in Lemma~\ref{lem:linearpack} is optimal. For \(\Phi(x,u)=xu\) and \(R=1\), the projected set with \(|\nabla_u\Phi|\le\varepsilon\le1\) is \([-\varepsilon,\varepsilon]\times[-1,1]\). Its \(t\)-packing number divided by \(t^{-1}\) tends to infinity whenever \(\varepsilon/t\to\infty\) and \(t\downarrow0\).
\end{example}

\section{The rank of a finite system of relations}\label{sec:rank}
Let \(S=(s_1,\ldots,s_n)\) and \(M=W^*(S)\) be as in Theorem~\ref{thm:main}. Fix a tuple \(Y=(y_1,\ldots,y_q)\in M_{\sa}^q\) and a column \(P=(P_1,\ldots,P_\ell)^{\mathsf T}\) of self-adjoint noncommutative polynomials in \(n+q\) formally self-adjoint letters. For a word \(w\) in the letters \(Z_1,\ldots,Z_{n+q}\), define its derivative in the \(j\)th letter by
\[
 \partial_j w=\sum_{w=uZ_jv}u\otimes v^{\op},
\]
where the sum is over the occurrences of \(Z_j\) in \(w\). Extend the derivative \(\partial_j\) by complex linearity. The evaluated Jacobian \(J\) and its blocks \(A,B\) are
\begin{equation}\label{eq:jacobian}
 J=\big((\partial_jP_i)(S,Y)\big)_{i,j}=[\,A\ B\,],\qquad
 A\in M_{\ell,n}(N),\quad B\in M_{\ell,q}(N).
\end{equation}
For \(j\ge1\), let \(\tau_j=\Tr_j\otimes\tau_N\) on \(M_j(N)\), where \(\Tr_j\) is the unnormalized matrix trace. For \(T\in M_{\ell,h}(N)\), write \(p_T=\Id_{(0,\infty)}(TT^*)\) for the range projection and \(\rk_N T=\tau_\ell(p_T)\) for the rank. Let \(\mu_T\) be the spectral distribution of \((TT^*)^{1/2}\) with respect to \(\tau_\ell\). Thus \(\mu_T\) has mass \(\ell\) and \(\mu_T((0,\infty))=\rk_N T\).

\begin{proposition}\label{prop:rank}
If \(P(S,Y)=0\), then the matrices \(J\) and \(B\) in \eqref{eq:jacobian} satisfy \(\rk_N J=\rk_N B\).
\end{proposition}

We first describe the matrix approximations used in the proof. For \(k\ge1\), let \(H_k=M_k(\C)_{\sa}\) with its real inner product \(\langle x,y\rangle_k=\tr_k(xy)\), where \(\tr_k=k^{-1}\Tr_k\). Write \(D=k^2=\dim_{\R}H_k\) and use \(\|\cdot\|_2\) for the induced norm on each direct sum \(H_k^j\). For \(x\in H_k^n\), \(y\in H_k^q\), the ordinary real derivative of \(P(x,y)\) is a map \(J_k=[\,A_k\ B_k\,]:H_k^{n+q}\to H_k^\ell\). Its complexification on \(M_k(\C)^{n+q}\) is the Jacobian \eqref{eq:jacobian} at \((x,y)\) acting by left and right multiplication. The complexification has exactly the same singular values, with the same multiplicities. In particular,
\begin{equation}\label{eq:super}
 D^{-1}\Tr_{\mathrm{End}(M_k(\C))}(h\mapsto ahb)
     =\tr_k(a)\tr_k(b),\qquad a,b\in M_k(\C).
\end{equation}
Indeed, summing the diagonal coefficients in the matrix-unit basis gives \eqref{eq:super}. Consequently the normalized trace \(D^{-1}\Tr((J_kJ_k^*)^v)\), for each integer \(v\ge0\), is a finite sum of products of normalized word traces in \((x,y)\). The same calculation applies to \(B_k\) and to mixed products of the blocks and their adjoints of compatible sizes.

Let \(X^{(k)}\in H_k^n\) consist of \(n\) independent Gaussian unitary ensemble (GUE) matrices normalized by \(\mathbb E\tr_k((X^{(k)}_i)^2)=1\). For each finite list of word traces of \((S,Y)\), each error \(\varepsilon>0\), and all sufficiently large \(k\), there is an event of probability at least \(1/2\) on which \(x=X^{(k)}\) admits \(y\in H_k^q\) approximating the specified traces within \(\varepsilon\). All the entries of \(x,y\) may be taken of operator norm at most a constant \(K\) depending only on the norms of the entries of \(S,Y\). To see this, use Kaplansky density in \(M=W^*(S)\) and norm density of polynomials in \(C^*(S)\) to choose self-adjoint polynomial tuples \(h_j(S)\) converging to \(Y\) in \(L^2(M)^q\) on a common operator-norm ball. Word traces converge under this approximation by telescoping each word and applying Cauchy--Schwarz. Fix \(j\) after fixing the list and the precision. The GUE moment theorem of Voiculescu \cite{Voi}, in its almost sure form \cite[(1.2)]{HT}, and norm convergence for each fixed polynomial \cite[Theorem~A]{HT} show that \(y=h_j(X^{(k)})\) has the required joint traces and norm bounds with probability tending to one. Thus \(K\) can be fixed before choosing the moment precision. On this ball the matrix blocks \(A_k,B_k\) have uniformly bounded norms, since \(\|ahb\|_2\le\|a\|\,\|b\|\,\|h\|_2\). Uniform polynomial approximation and \eqref{eq:super} therefore transfer any fixed finite list of continuous spectral tests on \(J_kJ_k^*\) or \(B_kB_k^*\).

\begin{proof}[Proof of Proposition~\ref{prop:rank}]
Assume that \(r=\rk_N B\) and \(b=\rk_N J-r>0\). Choose \(\gamma>0\) so that \(\mu_J([2\gamma,\infty))>r+7b/8\). For each \(0<\alpha<\gamma/2\), the preceding approximation gives, on a set of base points \(x=X^{(k)}\) of probability at least \(1/2\), a tuple \(y\) satisfying
\begin{align}
 \|P(x,y)\|_2&\le\alpha, \nonumber\\
 D^{-1}\dim\Id_{[\gamma,\infty)}((J_kJ_k^*)^{1/2})&\ge r+3b/4,\label{eq:spec}\\
 D^{-1}\dim\Id_{(\alpha,\infty)}((B_kB_k^*)^{1/2})&\le r+b/8.\nonumber
\end{align}
Here and below \(\dim e\) denotes the real dimension of the range of a real orthogonal projection \(e\). For the lower bound in \eqref{eq:spec}, use a continuous function between \(0\) and \(1\) which vanishes on \([0,\gamma]\) and equals \(1\) on \([2\gamma,\infty)\). For the upper bound, use a continuous function between \(0\) and \(1\) which vanishes on \([0,\alpha/2]\) and equals \(1\) on \([\alpha,\infty)\). Its integral against \(\mu_B\) is at most \(r\). The residual bound follows from \(P(S,Y)=0\).

Let \(V_{x,y}\) be the intersection of the high spectral subspace for \(J_kJ_k^*\) in \eqref{eq:spec} and the orthogonal complement of the high spectral subspace for \(B_kB_k^*\) in \eqref{eq:spec}. The dimension estimate for intersections of subspaces gives \(\dim V_{x,y}\ge(b/2)D\). Every \(w\in V_{x,y}\) satisfies
\[
 \|B_k^*w\|_2\le\alpha\|w\|_2,\qquad
 \|A_k^*w\|_2^2=\|J_k^*w\|_2^2-\|B_k^*w\|_2^2
 \ge(\gamma^2-\alpha^2)\|w\|_2^2.
\]
Let \(c=\min(1,\gamma/2)\). The image under \(A_k^*\) of the unit ball of \(V_{x,y}\) contains the ball of radius \(c\) in the subspace \(A_k^*V_{x,y}\). Hence, for \(0<t<c\), the image contains at least \((c/t)^{bD/2}\) points separated by \(t\), by the volume bound for a covering by \(t\)-balls.

In orthonormal real coordinates on \(H_k^n\), the Gaussian vector \(X^{(k)}\) has covariance \(D^{-1}I_{nD}\). A ball of radius \(t\) therefore has probability at most \((C_n t)^{nD}\), with \(C_n=\max(1,\sqrt{e/n})\). Indeed the density is at most \((D/2\pi)^{nD/2}\) and \(\omega_j\le(2\pi e/j)^{j/2}\). A maximal \(t\)-separated set of good base points \(x\) has at least \(\frac12(C_nt)^{-nD}\) points. Above each base point choose the separated points \(z=A_k^*w\) just constructed. The resulting pairs \((x,z)\) are \(t\)-separated in \(H_k^n\oplus H_k^n\).

For the upper bound on the number of these pairs, apply Lemma~\ref{lem:linearpack} to the real polynomial
\begin{equation}\label{eq:potential}
 \Phi_k(x,y,w)=\sum_{i=1}^\ell\tr_k(w_iP_i(x,y)),\qquad
 (x,y,w)\in H_k^n\oplus H_k^q\oplus H_k^\ell .
\end{equation}
The trace in \eqref{eq:potential} is real because \(w_i\) and \(P_i(x,y)\) are self-adjoint. Its gradients in \(x\) and in \((y,w)\) are \(A_k^*w\) and \((B_k^*w,P(x,y))\), respectively. Set \(\alpha=t/2\). On the chosen pairs the auxiliary gradient has norm at most \(t\). The norms of \(x,(y,w),z\) are bounded in terms of \(P,K,n,q,\ell\) only. The degree of \(\Phi_k\) is at most \(1+\max_i\deg P_i\). Lemma~\ref{lem:linearpack} and the two lower bounds yield
\begin{equation}\label{eq:compare}
 \tfrac12(C_n t)^{-nD}(c/t)^{bD/2}
 \ \le\ C^{(n+q+\ell)D}t^{-nD}.
\end{equation}
First choose \(t>0\) so small that \(C_n^{-n}(c/t)^{b/2}>2C^{n+q+\ell}\). For this fixed \(t\), choose \(k\) large enough for \eqref{eq:spec}. The two sides of \eqref{eq:compare} are then incompatible. This proves \(\rk_N J=\rk_N B\).
\end{proof}

\section{Extension to the von Neumann algebra}\label{sec:extension}
\begin{proof}[Proof of Theorem~\ref{thm:main}]
For each exact system \(P(S,Y)=0\), the range projections of the two matrices \(B\) and \(J=[\,A\ B\,]\) satisfy \(p_B\le p_J\). Proposition~\ref{prop:rank} and faithfulness of \(\tau_\ell\) yield \(p_B=p_J\). Let \(B^\dagger\in M_{q,\ell}(E)\) denote the Moore--Penrose inverse of \(B\). Since \(BB^\dagger=p_B\), the equality \(p_BA=A\) gives
\begin{equation}\label{eq:finite-sol}
 A\eta+B\zeta=0,\qquad \zeta=-B^\dagger A\eta,\quad \eta\in E^n .
\end{equation}
In particular, each finite system of differentiated relations admits a solution with the prescribed values \(\eta\) on \(S\). We now pass from \eqref{eq:finite-sol} to a derivation on the entire algebra \(M\).

Let \(\mathfrak A=\C\langle Z_h:h\in M_{\sa}\rangle\) be the free algebra on a distinct formally self-adjoint letter \(Z_h\) for each \(h\in M_{\sa}\). The evaluation map \(\pi:\mathfrak A\to M\), \(\pi(Z_h)=h\), is onto since \(a=\Re a+i\Im a\) for \(a\in M\). Let \(I=\ker\pi\), let \(G=\bigoplus_{h\in M_{\sa}} Ee_h\) be the free left \(E\)-module on the symbols \(e_h\), and let
\[
 K=G\Big/\sum_{P=P^*\in I} E\,r_P,\qquad
 r_P=\sum_{h\in M_{\sa}}(\partial_hP)(\pi)e_h .
\]
Every relation row \(r_P\) has finite support. Write \(d h\) for the class of \(e_h\) in \(K\). A left \(E\)-linear map \(\varphi:K\to E\) gives a derivation on \(\mathfrak A\) by assigning \(\varphi(dh)\) to \(Z_h\) and applying \eqref{eq:leibniz} through \(\pi\). For a self-adjoint polynomial \(P\in I\), its derivative is \(\sum_h(\partial_h P)(\pi)\varphi(dh)=0\), because the row \(r_P\) has class zero in \(K\). Since \(I\) is a \(*\)-ideal, the real and imaginary parts of each member of \(I\) also belong to \(I\). Thus the derivation descends to \(\mathfrak A/I=M\). Conversely, the values of every \(\delta\in\Der_\C(M,E)\) define a left \(E\)-linear map \(\varphi:K\to E\).

We claim that \(ds_1,\ldots,ds_n\) are left \(E\)-linearly independent in \(K\). Suppose \(\sum_i b_i\,ds_i=0\). In the free module \(G\) the row \(\sum_i b_i e_{s_i}\) is a finite left \(E\)-linear combination of relation rows \(r_{P_1},\ldots,r_{P_\ell}\). List as \(Y\) all the letters occurring in these rows other than \(s_1,\ldots,s_n\). Fix \(\eta\in E^n\). Equation~\eqref{eq:finite-sol} supplies a column \(\zeta\in E^q\) for this list. Evaluate the basis of \(G\) at \(\eta\) on \(S\), at \(\zeta\) on \(Y\), and at zero off \(S\cup Y\). The resulting functional annihilates each of the \(\ell\) specified rows. Therefore \(\sum_i b_i\eta_i=0\). Taking \(\eta\) to be each standard basis vector proves \(b_1=\cdots=b_n=0\).

The submodule \(\sum_i E\,ds_i\subset K\) is consequently free, and \(ds_i\mapsto\eta_i\) defines a left \(E\)-linear map on this submodule. The ring \(E=\Aff(N)\) is self-injective as a left module over itself \cite[Theorem~3.1]{Thom}. Extend the map to a left \(E\)-linear map \(\varphi:K\to E\). The corresponding derivation \(\delta:M\to E\) satisfies \(\delta(s_i)=\eta_i\) for every \(i\).
\end{proof}

\begin{proof}[Proof of Corollary~\ref{cor:main}]
Set \(D_a=L_a-R_a\) for \(a\in M\). Apply Theorem~\ref{thm:main} with \(\delta(s_1)=0\) and \(\delta(s_2)=1_E\). In the product spectral representation of \(W^*(s_1)\bar\otimes W^*(s_1)^{\op}\), the operator \(D_{s_1}\) is multiplication by \(\lambda-\mu\). The semicircle distribution has no atoms, so \(\lambda-\mu\ne0\) almost everywhere for the product measure. Hence \(D_{s_1}\) is invertible in \(E\). If \(\delta=\delta_\xi\) were inner, \(D_{s_1}\xi=0\) would give \(\xi=0\), contrary to \(\delta(s_2)=1_E\). This proves \(H^1_{\alg}(M,E)\ne0\). The group factor \(M=L(\F_n)\) has separable predual, so Theorem~\ref{dim:dichotomy} below gives \(\beta_1^{(2)}(M,\tau)=\infty\).
\end{proof}

Surjectivity of \(\ev_S\) also gives the bound \(\beta_1^{(2)}(L(\F_n))\ge n-1\) directly. Indeed, the image of \(\Inn(M,E)\) under \(\ev_S\) is the image of \(\xi\mapsto (D_{s_i}\xi)_{i=1}^n\), which is a copy of \(E\) because \(D_{s_1}\) is invertible. Thus \(H^1_{\alg}(M,E)\) maps onto a module of dimension \(n-1\). The outer derivation constructed above is discontinuous from the operator norm on \(M\) to the measure topology on \(E\), by \cite[Theorem~1]{PV}.

The preceding lower bound gives another proof of Corollary~\ref{cor:main}.
Shlyakhtenko's isomorphism theorem for semicircular \(C^*\)-algebras
\cite[Corollary~5.3]{Shlyakhtenko2026} yields a trace-preserving
isomorphism \(L(\F_m)\cong L(\F_2)\) for every integer \(m\ge2\). Hence
\[
 \beta_1^{(2)}(L(\F_2))=\beta_1^{(2)}(L(\F_m))\ge m-1
 \qquad(m\ge2).
\]
Letting \(m\to\infty\) shows that the common Betti number is infinite.
The nonvanishing of the first algebraic cohomology group follows from
\eqref{eq:duality}.

\subsection{A cycle with bounded coefficients}
There is an explicit representative on the homology side of \eqref{eq:duality}. Let \(n=2\). In the ring \(E\) define
\[
 T=D_{s_2}D_{s_1}^{-1},\qquad h=(1+TT^*)^{-1}.
\]
The operator \(h\in N\) has full support, and \(hT\in N\) because \(\sup_{t\ge0} t/(1+t^2)=1/2\). In the bar complex of \(M\) with coefficients in \(N\) \cite[Section~2.3]{CS}, consider the chain
\begin{equation}\label{eq:cycle}
 z=h\otimes s_2-hT\otimes s_1\in N\otimes_{\C}M .
\end{equation}
The first boundary \(b_1(c\otimes a)=cD_a\) gives \(b_1(z)=0\). For a derivation \(\delta\) as in the proof of Corollary~\ref{cor:main}, pair a chain \(c\otimes a\) with \(\delta\) by \(c\,\delta(a)\). The second boundary is given by
\[
 b_2(c\otimes a\otimes b)
   =cL_a\otimes b-c\otimes ab+cR_b\otimes a .
\]
The pairing annihilates \(\operatorname{im}b_2\) by the Leibniz identity \eqref{eq:leibniz}, whereas the value of the pairing on \(z\) is \(h\ne0\). In fact \(c\mapsto[cz]\) embeds the left module \(N\) into the first homology group, since the vanishing of the class \([cz]\) would give \(ch=0\), and \(h\) has full support.

\subsection{A spectral projection}
An elementary relation explains why affiliated coefficients occur in the extension. Prescribe the value \(\delta(s_1)=1_E\) and put \(e=\Id_{[0,\infty)}(s_1)\). Differentiating \(s_1e=es_1\) gives \(D_{s_1}\delta(e)=D_e\), so every extension has
\[
 \delta(e)(\lambda,\mu)=
 \frac{\Id_{[0,\infty)}(\lambda)-\Id_{[0,\infty)}(\mu)}{\lambda-\mu}
\]
in the same product spectral representation. Since the semicircle density is bounded below near zero, the square integral of \(\delta(e)\) dominates a positive constant times \(\int_0^\epsilon\int_0^\epsilon(u+v)^{-2}\,du\,dv=\infty\). Thus the forced value \(\delta(e)\) belongs to \(E\) and lies outside \(L^2(N)\).

\section{The dimension of the first cohomology group}\label{dim:sec}

\begin{theorem}\label{dim:dichotomy}
Let $(M,\tau)$ be a $\mathrm{II}_1$ factor with separable predual, let
$N=M\bar\otimes M^{\op}$, and let $E=\Aff(N)$.
If $\beta_1^{(2)}(M)<\infty$, then every algebraic derivation
$\delta\colon M\to E$ is inner. Consequently, the existence of an outer
derivation $\delta\colon M\to E$ implies $\beta_1^{(2)}(M)=\infty$.
\end{theorem}

\begin{proof}
Write $\tau_N=\tau\otimes\tau^{\op}$ and
$\mathcal D=\Der_{\C}(M,E)$, with right $E$-action
$(\delta b)(x)=\delta(x)b$. By Thom's duality
\cite[Corollaries~3.4 and~3.7]{Thom},
$\dim_N(\mathcal D/\Inn(M,E))=\beta_1^{(2)}(M)$.
The map from $E$ onto $\Inn(M,E)$ taking an implementer to its inner
derivation is right $E$-linear. Hence
$\dim_N\mathcal D\leq 1+\beta_1^{(2)}(M)<\infty$, by additivity
of the extended von Neumann dimension \cite[Chapter~6]{Luck}.

For a finite subset $F\subset M$, let $V_F$ be the image of the
evaluation map $\mathcal D\to E^F$. Put
$r=\sup_F\dim_N V_F$ and choose increasing finite sets $F_j\subset M$
with $\dim_N V_{F_j}\to r$. The number $r$ is finite since
$\dim_N V_F\leq\dim_N\mathcal D$. Put $S=\bigcup_j F_j$.
For $x\in M\setminus S$, the module of values
$\{\delta(x):\delta\in\mathcal D,\ \delta|_S=0\}$ embeds into the
kernel of the restriction $V_{F_j\cup\{x\}}\to V_{F_j}$ for every $j$.
The restriction is onto, so its kernel has dimension at most
$r-\dim_N V_{F_j}$. Thus the module of values at $x$ has dimension zero.
Every nonzero right $E$-submodule of $E$ has positive dimension. Indeed,
if the submodule contains $b\ne0$, it contains
$(bb^\dagger)N$, where $b^\dagger$ is the Moore--Penrose inverse in $E$.
It follows that a derivation vanishing on $S$ vanishes on $M$.

We show that the values on the countable set $S$ give a measurable
description of each derivation. Let $K$ be the left $E$-module with
generators $dx$ for $x\in M$. For $x,y\in M$ and $\lambda\in\C$ the defining relations are
\[
 d(\lambda x+y)=\lambda\,dx+dy,\qquad
 d(xy)=L_x\,dy+R_y\,dx .
\]
The assignment $dx\mapsto\delta(x)$ identifies $\mathcal D$ with
$\operatorname{Hom}_E(K,E)$. Write $U=\sum_{s\in S} E\,ds$. For
$x\in M$, define the left ideal $I_x=\{b\in E:b\,dx\in U\}$.
Its right annihilator is zero. Indeed, if $I_x w=0$, the formula
$u+b\,dx\mapsto bw$ defines a left $E$-linear map $U+E\,dx\to E$.
Left self-injectivity of $E$ \cite[Theorem~3.1]{Thom} extends the map
to $K$. The corresponding derivation vanishes on $S$ and takes the
value $w$ at $x$, so $w=0$.

It follows that $I_x$ contains projections of trace arbitrarily close
to $1$. To see this, note that $b^\dagger b\in I_x$ for $b\in I_x$.
The complement of $\bigvee_{b\in I_x}b^\dagger b$ belongs to the
right annihilator of $I_x$ and hence is zero. Finite joins of projections
in $I_x$ belong to $I_x$: the join of $p_1,\ldots,p_m$ is $h^\dagger h$
for $h=\sum_i p_i$. Normality of $\tau_N$ gives the assertion.

Fix $\delta\in\mathcal D$. For $x\in M$ and $y\in E$, the equality
$y=\delta(x)$ holds if and only if, for every integer $k\ge1$, there
exist a projection $e\in N$, a finite list $s_1,\ldots,s_m\in S$,
and coefficients $a_1,\ldots,a_m\in E$ such that
\begin{equation}\label{dim:witness}
 \tau_N(e)>1-2^{-k},\qquad e\,dx=\sum_{i=1}^m a_i\,ds_i\ \text{in }K,
 \qquad ey=\sum_{i=1}^m a_i\delta(s_i).
\end{equation}
For $y=\delta(x)$, take $e\in I_x$ and evaluate the middle identity
at $\delta$. Conversely, the last two identities give $e(y-\delta(x))=0$.
The projection onto the closure of the range of $y-\delta(x)$ therefore
has trace at most $2^{-k}$ for every $k$, so $y=\delta(x)$.

Equip $M$ with the standard Borel structure given by the $L^2$ metric
on its closed operator-norm balls, and equip $E$ with the Borel structure
of its Polish measure topology. The bounded balls of $M$ are Polish in
the $L^2$ metric because $M$ has separable predual. Addition,
multiplication, and the two embeddings of $M$ into $N$ are Borel for
these structures. Each equality in $K$ in \eqref{dim:witness} can be witnessed by expressing the difference as a finite sum of left multiples
of the defining relations of $K$, grouping the terms with the same label
$dz$, and comparing the coefficients in $E$. For each list length,
equality of labels in $M$ and equality of coefficients in $E$ are Borel
tests. The values $\delta(s)$ for $s\in S$ form a fixed countable table.
Consequently the graph of $\delta$ is analytic because it is the projection of
the Borel set of pairs $(x,y)$ together with a sequence of witnesses
for \eqref{dim:witness}, one for each $k$. For any Borel set $B\subset E$,
both $\delta^{-1}(B)$ and $M\setminus\delta^{-1}(B)$ are projections
of analytic sets, obtained by restricting this graph to $M\times B$
and $M\times(E\setminus B)$, respectively. Souslin's theorem
\cite[Theorem~14.11]{Kechris} implies that $\delta^{-1}(B)$ is Borel.

Every set in the indicated Borel structure on $M$ is norm Borel. Indeed,
on each closed norm ball the norm topology is finer than the $L^2$
topology, and countably many such balls cover $M$. Thus $\delta$ is a
Baire measurable homomorphism from the additive group $(M,\|\cdot\|)$
to the additive group $E$ in the measure topology. The domain is a
Baire group by completeness of the operator norm, and the target is
separable. Pettis's automatic continuity theorem
\cite[Theorem~9.10]{Kechris} gives norm-to-measure continuity of $\delta$.
The vanishing theorem of Popa and Vaes \cite[Theorem~1]{PV} now implies
that $\delta$ is inner.
\end{proof}

\section{An estimate in terms of free entropy dimension}\label{sec:entropy}
Let \(X=(x_1,\ldots,x_n)\) be a bounded self-adjoint tuple in a tracial von Neumann algebra, and now let \(M=W^*(X)\). Assume that \(X\) admits matricial microstates. Denote by \(\delta_0(X)\) the modified microstates free entropy dimension and by \(\ev_X(\Der_\C(M,E))\) the right \(E\)-module of value columns at \(X\).

\begin{proposition}\label{prop:entropy}
The dimension of the module of value columns satisfies
\[
 \delta_0(X)\le \dim_N\ev_X(\Der_\C(M,E)).
\]
If \(M\) is a \(\mathrm{II}_1\) factor and \(\delta_0(X)>1\), then \(\beta_1^{(2)}(M,\tau)=\infty\).
\end{proposition}

\begin{proof}
Fix a finite tuple \(Y\in M_{\sa}^q\) and a finite system of actual relations \(P(X,Y)=0\). Write \(J=[\,A\ B\,]\) for the evaluated Jacobian at \((X,Y)\) and \(b=\rk_N J-\rk_N B\). We first show that \(\delta_0(X)\le n-b\). This holds for \(b=0\) by the ambient dimension bound. If \(b>0\), fix \(0<b'<b\). At a given matrix size \(k\), take a joint microstate for \((X,Y)\) with an operator-norm cutoff greater than the norms of the entries of \((X,Y)\). Such a microstate approximates a finite list of word traces to a prescribed tolerance. For each fixed \(t>0\) choose the word list and the tolerance to enforce the residual and spectral tests in the proof of Proposition~\ref{prop:rank}. The two spectral cutoffs may be chosen with total loss less than \(b-b'\). The high spectral cutoff \(\gamma>0\) depends on \(b'\) and \(J\) only. It follows that each projected microstate in the \(X\)-coordinates admits a ball of \(A_k^*w\)-values of radius \(c>0\) and dimension at least \(b'k^2\), as in that proof.

Let \(\Gamma_X\) be the projection onto the \(X\)-coordinates of the resulting joint microstate set. Form a \(t\)-separated subset above each point of a maximal \(t\)-separated subset of \(\Gamma_X\). Applying Lemma~\ref{lem:linearpack} to the polynomial \(\Phi_k\) in \eqref{eq:potential} gives
\[
  \Pack_t(\Gamma_X)(c/t)^{b'k^2}
    \le C^{(n+q+\ell)k^2}t^{-nk^2},
\]
with \(C,c\) independent of \(t,k\) for \(t\) sufficiently small. Take \(k^{-2}\) times the logarithm, then the upper limit in \(k\), and then the infimum over the moment conditions. Divide by \(|\log t|\) and take the upper limit as \(t\downarrow0\). Jung's packing formula, including the presence of \(Y\subset W^*(X)\), identifies the resulting left packing exponent with \(\delta_0(X)\) \cite[Remark~2.3 and Corollary~2.4]{Jung}. Changing the packing scale by a fixed constant has no effect on this exponent. We obtain \(\delta_0(X)\le n-b'\) and hence \(\delta_0(X)\le n-b\).

For the given finite system let \(p_P\in M_n(N)\) be the kernel projection of \((1-p_B)A\). Solvability of \(A\eta+B\zeta=0\) over \(E\) is equivalent to \(\eta\in p_PE^n\), since \(BB^\dagger=p_B\). The operator \((1-p_B)A\) has closed range span \((p_J-p_B)L^2(N)^\ell\), so \(\tau_n(p_P)=n-b\). On combining any two finite systems and taking the union of their auxiliary letters, one obtains a system whose projection is below both original projections. For \(p=\bigwedge_P p_P\), normality of the trace therefore gives \(\tau_n(p)\ge\delta_0(X)\).

The module \(\ev_X(\Der_\C(M,E))\) is exactly \(pE^n\). Indeed, the values of any derivation satisfy every finite system. Conversely, any \(\eta\in pE^n\) annihilates every left \(E\)-linear relation among \(dx_1,\ldots,dx_n\) in the module \(K\) of Section~\ref{sec:extension}. Each such relation uses finitely many polynomial rows, and the values \(\eta\) extend across those rows by the definition of \(p\). Left self-injectivity of \(E\) extends the corresponding functional to \(K\). This proves the desired dimension estimate. Repeated entries of \(X\) can equivalently be kept as separate formal letters with the equality relations between their copies.

If \(\delta_0(X)>1\), the module of values has dimension greater than one, while the values of inner derivations have dimension at most one. Thus \(M\) has an outer derivation. The algebra \(M=W^*(X)\) has separable predual, and Theorem~\ref{dim:dichotomy} gives the last assertion.
\end{proof}

\appendix

\section{Logarithmic integrability for a finite system}
\label{q:section}

Let \(S=(s_1,\ldots,s_n)\) be a standard free semicircular system, let
\(M=W^*(S)\), and write \(N=M\bar\otimes M^{\op}\) and \(E=\Aff(N)\).
For a self-adjoint tuple \(Y=(y_1,\ldots,y_q)\) in \(M\) and a tuple
\(P=(P_1,\ldots,P_\ell)\) of self-adjoint noncommutative polynomials,
denote by \(J=[A\ B]\) the Jacobian \eqref{eq:jacobian} of \(P\)
evaluated at \((S,Y)\).
The columns of \(A\) correspond to \(S\), and the columns of \(B\) to \(Y\).
Put \(\tau_m=\Tr_m\otimes\tau_N\), where \(\Tr_m\) is the unnormalized
trace. We use the sum of squares of the component \(L^2\)-norms to
define \(\|P(S,Y)\|_2^2\).

\begin{theorem}\label{q:bound}
Fix the polynomial tuple \(P\), the integers \(n,q,\ell\), and \(K>0\).
There exists \(C<\infty\) such that, whenever \(\max_j\|y_j\|\le K\),
\(0<t<1\), and \(\|P(S,Y)\|_2\le t/2\), one has
\begin{equation}\label{q:trace}
 \tau_n\!\left(\log\!\left(1_n+A^*(BB^*+t^2 1_\ell)^{-1}A\right)\right)
       \le C .
\end{equation}
If \(P(S,Y)=0\), then \(C_*=B^\dagger A\in M_{q,n}(E)\) satisfies
\begin{equation}\label{q:solution}
 BC_*=A,\qquad \tau_n\!\left(\log(1_n+C_*^*C_*)\right)\le C.
\end{equation}
Here \(B^\dagger\) denotes the Moore--Penrose inverse of \(B\).
\end{theorem}

\begin{proof}
Write \(L_t\) for the left-hand side of \eqref{q:trace}. Give
\(H_k=M_k(\C)_{\sa}\) the real inner product
\(\langle x,y\rangle=k^{-1}\Tr_k(xy)\), and put \(D=k^2\). Take \(x\)
to be an independent GUE \(n\)-tuple normalized to covariance \(D^{-1}I\)
in orthonormal coordinates. For fixed \(t,\varepsilon>0\) as above and
sufficiently large \(k\), a set of such tuples \(x\) of probability at
least \(1/2\) admits tuples \(y\in H_k^q\) with a common operator-norm
bound, \(\|P(x,y)\|_2\le t\), and
\begin{equation}\label{q:transfer}
 \frac1D\Tr\log(1_{nD}+t^{-2}F_kF_k^*)\ge L_t-\varepsilon,
 \qquad F_k=A_k^*(1_{\ell D}+t^{-2}B_kB_k^*)^{-1/2}.
\end{equation}
The matrices \(A_k,B_k\) are the real derivative blocks of \(P\) at
\((x,y)\). To see \eqref{q:transfer}, approximate \(Y\) in \(L^2\) by
self-adjoint polynomial tuples \(h(S)\) in a fixed norm ball, using
Kaplansky density. At fixed \(t\), approximate
\(s\mapsto(1+t^{-2}s)^{-1/2}\) uniformly on the spectral interval
for \(B_kB_k^*\), then approximate the logarithm on the resulting
compact positive interval. The normalized endomorphism trace of
\(H\mapsto aHb\) on \(M_k(\C)\) is \(\tr_k(a)\tr_k(b)\).
Consequently each polynomial test is a sum of products of word traces.
Choose the approximation \(h\) first and put \(y=h(x)\).
GUE moment and norm convergence \cite{HT} give the two inequalities
with probability tending to one. The norm bound for \((x,y)\) depends
only on \(n,q,K\).

For each such \(x\), let \(w=(1_{\ell D}+t^{-2}B_kB_k^*)^{-1/2}v\),
where \(|v|\le1\). Then \(|w|\le1\), \(|B_k^*w|\le t\), and
\(z=A_k^*w\) ranges over the ellipsoid \(F_k\{v:|v|\le1\}\).
Let \(\sigma_1,\ldots,\sigma_{nD}\) be the axis lengths, with zeros
included. Projecting onto the axes of length \(>t\) and comparing
volumes with a ball of radius \(t\) gives a \(t\)-packing of size at least
\[
 \prod_i\max(1,\sigma_i/t)
 \ \ge\ 2^{-nD/2}\det(1_{nD}+t^{-2}F_kF_k^*)^{1/2}.
\]
The probability of any base ball of radius \(t\) is at most
\((C_nt)^{nD}\), by the maximum Gaussian density times the ball volume.
Thus at least \(\tfrac12(C_nt)^{-nD}\) mutually \(t\)-separated
base tuples \(x\) are available. Combine the packing at each base tuple
with \eqref{q:transfer}.

An upper bound on the resulting packing of pairs \((x,z)\) follows
from Lemma~\ref{lem:linearpack}, applied to the real polynomial
\(\Phi(x,y,w)=\sum_i\tr_k(w_iP_i(x,y))\). Its main gradient equals
\(z=A_k^*w\). Its auxiliary gradient \((B_k^*w,P(x,y))\) has norm
at most \(\sqrt2\,t\). All three norms \(|x|,|(y,w)|,|z|\) have fixed
bounds, and \(\deg\Phi\le 1+\max_i\deg P_i\). Comparing the upper and
lower bounds gives
\[
 \tfrac12(C_nt)^{-nD}
 e^{D(L_t-\varepsilon-n\log 2)/2}
 \ \le\ C_0^{(n+q+\ell)D}t^{-nD}.
\]
Taking \(k\to\infty\) and then \(\varepsilon\downarrow0\) proves
\eqref{q:trace}, with a constant independent of \(Y,t\).

Suppose \(P(S,Y)=0\) and set \(p=1_{\{0\}}(BB^*)\).
The positive operator \(K_t=A^*(BB^*+t^2 1_\ell)^{-1}A\) dominates
\(t^{-2}A^*pA\). If \(A^*pA\ne0\), a spectral projection of \(A^*pA\)
on some interval \([c,\infty)\), \(c>0\), has positive trace. Operator
monotonicity of the logarithm then contradicts \eqref{q:trace} as
\(t\downarrow0\). Hence \(pA=0\) and \(BC_*=BB^\dagger A=A\).
Put \(C_t=B^*(BB^*+t^2 1_\ell)^{-1}A\). One has \(C_t^*C_t\le K_t\)
and \(C_t=f_t(B^*B)C_*\), where \(f_t(s)=s/(s+t^2)\).
The multipliers \(f_t(B^*B)\) converge in measure to the support of
\(B^*B\), which fixes \(C_*\). Therefore \(C_t\to C_*\) in measure.
Fatou's lemma for generalized singular values \cite{FackKosaki} and
\eqref{q:trace} prove \eqref{q:solution}.
\end{proof}

\subsection{Optimality of the logarithmic exponent}
The choice of the finite system \(P\) matters in Theorem~\ref{q:bound}.
The next example fixes a single system for all exponents \(1+\varepsilon>1\).

\begin{proposition}\label{q:sharp}
Let \(M=W^*(s_0,s_1)\) for a standard free semicircular pair \((s_0,s_1)\).
There exist a finite tuple \(Y\in M_{\sa}^q\), a self-adjoint polynomial
system \(P(s_0,s_1,Y)=0\), and an index \(m\le q\) with the following
property. Every solution \(\eta\in E^q\) of \(A(0,1)^{\mathsf T}+B\eta=0\)
satisfies
\[
    \tau_N\!\left((\log(1+|\eta_m|))^{1+\varepsilon}\right)=\infty
             \qquad(\varepsilon>0).
\]
\end{proposition}

\begin{proof}
Write \(L_x=x\otimes1,\ R_x=1\otimes x^{\op}\), and \(D_x=L_x-R_x\).
A dot denotes the entry assigned to a variable by a solution of the
linearized equations. Each complex variable stands for two self-adjoint
variables. Each complex identity stands for its self-adjoint real and
imaginary parts. In particular, \(\dot{(z^*)}\) denotes the derivative
assigned to \(z^*\) by this convention.

\smallskip\noindent\emph{A countable partition in a finite system.}
For \(j\ge1\) put
\[
 d=j+1,\qquad R=2^d,\qquad h=R^d-1,\qquad
 n_j=2^{d^2+10d+10},\qquad a_j=\frac{c}{d^2n_j},
\]
where \(c\) is chosen so that \(\sum_j n_j a_j=1\).
The parameters \(d,R,h\) in this construction refer to the current block \(j\).
Partition \(1\) in \(A_0=W^*(s_0)\) into projections
\(e_{j,l}\) of trace \(a_j\), \(0\le l<n_j\).
Let \(e_j=e_{j,0}\), \(p=\sum_j e_j\), and
\(\Pi=\sum_j L_{e_j}R_{e_j}\).
We first prescribe \((\dot s_0,\dot s_1)=(0,\Pi)\).
For any of the finitely many coefficients \(b\in A_0\) used below, add
the relation \([s_0,b]=0\). Its derivative \(D_{s_0}\dot b=0\)
forces \(\dot b=0\), since \(s_0\) has diffuse law. We call such
a coefficient \(b\) fixed.

Enumerate the projections \(e_{j,l}\) as \(e_i\), \(i\ge1\). Let
\((f_i)_i\) be a permutation of \((e_i)_i\). Use \(\lambda_i=i^{-1}\)
and \(b_i=2^{-i}\) as eigenvalues of diagonal coefficients
\(h_L,b_L\) on \(e_i\) and \(h_R,b_R\) on \(f_i\), respectively.
Add the two identities
\begin{equation}\label{q:read}
 h_L Z=Z h_R,\qquad b_L s_1 b_R=Z+h_L T-T h_R .
\end{equation}
Actual bounded values are
\[
 Z=\sum_i b_i^2e_i s_1 f_i,\qquad
 T=\sum_{i\ne k}\frac{b_i b_k}{\lambda_i-\lambda_k}e_i s_1 f_k .
\]
The series for \(T\) is norm convergent because
\(|\lambda_i-\lambda_k|^{-1}\le ik\). The kernel projection of
\(L_{h_L}-R_{h_R}\) equals \(k_0=\sum_i L_{e_i}R_{f_i}\).
Applying \(k_0\) to the linearization of \eqref{q:read} gives
\(\dot Z=k_0L_{b_L}R_{b_R}\Pi\).

First, let \(f_i\) precede \(e_i\) in the cyclic ordering within its
block \(j\). Then \(\dot Z=\dot{(Z^*)}=0\). The polar unitary \(U\) of
\(Z\) carries each projection \(f_i\) onto \(e_i\). Indeed, put
\(r=e_i+f_i\) and \(s'=rs_1r/\sqrt{\tau(r)}\). In \(rMr\), the
element \(s'\) is standard semicircular and free from \(rA_0r\) for the
normalized trace, by free compression \cite[\S\S 1.11--1.14]{NS}.
The centered symmetry \(b=e_i-f_i\) has the property that \(s'\) and
\(bs'b\) are free, as follows by expanding alternating centered words.
Thus \(s'-bs'b=2(e_i s'f_i+f_i s'e_i)\) has diffuse semicircular law.
It follows that \(e_i s_1f_i\) has full source and range supports \(f_i,e_i\).
Impose \(H^2=Z^*Z,\ Z=UH,\ U^*U=UU^*=1\), with \(H=|Z|\).
The positivity and full support of \(H\) make both \(L_H+R_H\)
and \(R_H\) invertible in \(E\). These identities force
\(\dot H=\dot U=\dot{(U^*)}=0\).

Second, use \eqref{q:read} with \(f_i=e_i\) and new variables \(Z,T\).
Define \(y=\sum_{i=(j,l)}e_i s_1 e_i/\sqrt{a_j}\), and impose \(Z=by\)
with fixed \(b=\sum_{i=(j,l)}b_i^2\sqrt{a_j}\,e_i\).
Free compression gives the semicircular law in each corner and hence
for \(y\). Now \(\dot y=\sum_j a_j^{-1/2}L_{e_j}R_{e_j}\).
Take \(h_1=\exp(2\pi i F_{\rm sc}(y))\), where \(F_{\rm sc}\) is the
semicircular distribution function. The relation \([y,h_1]=0\)
gives \(\dot h_1=D_y^{-1}D_{h_1}\dot y\). Choose a fixed unitary
\(h_0\in A_0\) with Haar law on each projection \(e_i\).
In the corner \(e_jMe_j\), the unitaries \(u=e_jh_0e_j\) and
\(v=e_jh_1e_j\) are free Haar unitaries. Writing \(y_j=e_jye_j\),
the coefficient
\[
 \zeta_j=R_{v^*}D_{y_j}^{-1}D_v
\]
is invertible over the coarse square of \(e_jMe_j\). The legs in this
formula have unit \(e_j\). For its normalized trace the singular values
satisfy \(\mu_s(\zeta_j^{-1})\le C_1/s\), \(0<s<1\).
Indeed \(\|D_{y_j}\|\le4\), and the probability that two independent
Haar points are within distance \(\alpha\) is \(O(\alpha)\).

\smallskip\noindent\emph{Matrix coefficients with zero derivative.}
In block \(j\) put \(v_a=U^a e_j,\ 0\le a<n_j\), and use
\(v_av_b^*\) as matrix units. Define \(V=(1-p)U\) and
\(H_0=\sum_{j,a}v_a(e_jh_0e_j)v_a^*\). Impose
\(H_0=ph_0+VH_0V^*\). Its linearization gives
\(\dot H_0=L_{V^m}R_{(V^*)^m}\dot H_0\) for every \(m\).
The contraction \(V\) acts as a truncated shift of length \(n_j\)
in block \(j\). Hence \(V^m\to0\) in measure, and multiplication in
\(E\) gives \(\dot H_0=0\). The law of \(H_0\) is Haar. Every bounded
block sum \(K\) of scalar matrices in the units \(v_av_b^*\) commutes
with \(H_0\). Imposing the single relation \([H_0,K]=0\) therefore
forces \(\dot K=0\). This observation provides the finitely many
fixed matrix coefficients in the construction below.

\smallskip\noindent\emph{Encoding entrywise multiplication.}
In block \(j\) store a rectangular array \(X=(x_{ik})\) with entries
in \(e_jMe_j\) as the single variable
\[
 [X]=\sum_{j,i,k}\sigma_j v_i x_{ik}v_k^*,\qquad
 \sigma_j=e^{-\sqrt h}.
\]
In the formula for \([X]\), the notation \(v_i,x_{ik},v_k\) denotes the partial isometries and the entry belonging to the summand indexed by \(j\).
Let \(\Sigma\) be the fixed scalar matrix with value \(\sigma_j\) on
block \(j\). Impose the row and column support equation
\([X]=p_I[X]p_J\), and impose commutation with the fixed block diagonal
having value \(1/(j+1)\) on block \(j\). The derivative entry is
\(\dot x_{ik}=\sigma_j^{-1}L_{v_i^*}R_{v_k}\dot{[X]}\).
Inserting the site projections in each matrix product shows that
the raw product \(C=AB\), imposed by \([A][B]=\Sigma[C]\),
obeys the usual matrix product rule on these derivative entries.
Their left support is dominated by \(L_{e_j}R_{e_j}\).

We describe a fixed number of polynomial identities which impose
entrywise multiplication \(Z=X\odot Y\) on \(h\)-by-\((h+1)\) arrays
and impose its product rule on derivative entries. First, diagonal
copies \(W=\operatorname{diag}_{r=1}^D[X]\) are specified by
\[
   H_IW=WH_J,\qquad W=U_IWU_J^*,\qquad I_1^*WJ_1=[X].
\]
The fixed matrices \(H_I,H_J\) label consecutive copies by \(r/D\);
\(U_I,U_J\) cyclically permute the copies, and \(I_1,J_1\) include the
first copy. The first equation removes off-diagonal blocks, the second
identifies diagonal blocks, and the third specifies their entries.
The equations have the same effect on derivative entries.

For sites with integer coordinates \(x,y\in[0,2R]^d\), let \(M_b\)
retain exactly the entries where \(k(x,y)=\|x-y\|^2-b=0\).
The integer \(b\in[0,4dR^2]\) may vary between copies. The vectors
\[
 a_x=(\|x\|^2-b,1,-2x_1,\ldots,-2x_d),\quad
 b_y=(1,\|y\|^2,y_1,\ldots,y_d)
\]
have scalar product \(k(x,y)\) and norms at most \(c_j=10dR^2\).
On \(d+2\) diagonal copies, multiply from the left and right by the
fixed matrices with entries
\(\delta_{xx'}(a_x)_r/c_j\) and
\(\delta_{yy'}(b_y)_r/c_j\), respectively. The resulting linear
map \(\mathcal K\) multiplies entry \((x,y)\) by
\(\kappa_j k(x,y)\), with \(\kappa_j=c_j^{-2}\).
Thus \(X=M_b(F)\) is imposed by the two equations
\[
  \mathcal K[X]=0,\qquad \kappa([F]-[X])=\mathcal K[T],
\]
where \(\kappa\) is the block scalar with value \(\kappa_j\). Take
\(t_{xy}=f_{xy}/k(x,y)\) off the zero set and zero on the zero set.
The same two equations select \(\dot f_{xy}\) on the zero set.

Use column indices \(J=\{0,\ldots,R-1\}^d\), row indices \(I=J\)
with one vertex removed, and middle indices
\(m\in\{0,\ldots,2R-2\}^d\). The Euclidean coordinates of the three
types of sites are \(2i,m,2k\). Take diagonal copies for the integers
\(0\le l\le d(R-1)^2\). In copy \(l\) use the fixed matrices
\[
   (B_R)_{km}=1_{\{\|2k-m\|^2=l\}},\qquad
   (B_L)_{mi}=1_{\{\|m-2i\|^2=l\}} .
\]
Use the bounded stored matrices \([B_R]\) and \([B_L]\) in the product equations \([A][B]=\Sigma[C]\). Form, in each copy, the arrays
\begin{align*}
 W&=M_{4l}(X), & T&=M_{4l}(Y),\\
 U'&=M_l(WB_R), & V'&=M_l(B_LT), & H_l&=M_{4l}(U'V').
\end{align*}
For each nonzero summand the endpoints have squared distance \(4l\)
and the middle point has squared distance \(l\) from both endpoints.
Equality in the triangle inequality forces \(m=i+k\). In particular,
the retained entries of \(U'\) and \(V'\) each have only the required
input index, respectively \(k=m-i\) and \(i=m-k\). Therefore
\((H_l)_{ik}=x_{ik}y_{ik}\) when \(l=\|i-k\|^2\), and is zero
otherwise. Summing diagonal copies uses the fixed row and column
summation matrices with coefficients \(D_l^{-1/2}\), where
\(D_l=1+d(R-1)^2\), and the fixed block scalar \(D_l^{-1}\) on the
output. This gives \(Z=X\odot Y\) and
\(\dot z_{ik}=L_{x_{ik}}\dot y_{ik}+R_{y_{ik}}\dot x_{ik}\).

These operations require at most
\(3(2R)^d[1+(1+dR^2)(d+3)]<n_j\) sites in each block.
The number of operations is independent of \(j\). The variable copy
counts occur only in the fixed matrices. For input entries of norm
at most one, all raw intermediate entries are bounded by a fixed
polynomial in \(n_j\), since each ordinary product sums at most \(n_j\)
terms, and division by a nonzero integer \(k(x,y)\) costs at most one.
The norm of an array is at most \(n_j\) times its maximum entry norm.
Since \(\log n_j=O(\log h)\), all stored variables and all scaled
scalar matrices define bounded elements of \(M\).

\smallskip\noindent\emph{The word in the last cell.}
Order the \(h\) rows and \(h+1\) columns cyclically from \(0\).
Let \(\mathcal S\) increase both indices by one, and let \(Q\) select
cell \((0,0)\), using fixed permutations and projections. The orbit
of \((0,0)\) has length \(H_j=h(h+1)\).
In root coordinates let \(U_c,U_c^{-}\) have respective entries
\(u,u^*\) at every cell. Their stored variables commute with \(H_0\),
so imposing commutation forces their derivatives to vanish. Let \(V_c\)
have \(v\) at \((0,0)\) and zero elsewhere. The stored seed \([V_c]\)
equals \(\Sigma K_1h_1K_2\), where the fixed partial isometries
\(K_1,K_2\) move root indices to the designated row and column.
Impose
\[
 G(X)=X\odot U_c\odot X\odot X\odot U_c^{-},\qquad
 [X]=[V_c]+(1-Q)\mathcal S[G(X)] .
\]
The actual entry at step \(r\) is the unitary \(w_r\), where
\(w_0=v\) and \(w_{r+1}=w_r u w_r^2u^*\).
Deleting the return to the initial cell gives an ordinary finite
induction on \(0\le r<H_j\) for the derivative entries as well.
Write \(\pi(w)=L_wR_{w^*}\) in the coarse root corner and set
\[
 B_0=1,\qquad B_{r+1}=
  \big(1+\pi(w_ru)+\pi(w_ruw_r)\big)B_r .
\]
The identity \(R_{(xy)^*}\dot{(xy)}
 =R_{x^*}\dot x+\pi(x)R_{y^*}\dot y\) and
\(R_{v^*}\dot v=a_j^{-1/2}\zeta_j\) show that the terminal derivative
entry of \([X]\), transported to the root, is
\begin{equation}\label{q:terminal}
  \sigma_j a_j^{-1/2}R_{w_{H_j-1}}B_{H_j-1}\zeta_j .
\end{equation}
The transport consists of partial isometries between equal-trace
corners and preserves the nonzero singular values under the ambient trace.

The exponent sum of \(v\) in \(w_r\) is \(3^r\), so \(w_r\) does
not commute with \(u\). The elements \(u,w_r\) form a free basis
of their two-generator free subgroup. The pair \(w_ru,w_ruw_r\)
is also a basis, as \(w_r=(w_ru)^{-1}(w_ruw_r)\).
The representation \(\pi\) preserves the free group trace on this
subgroup because the normalized coarse corner trace of \(\pi(g)\) equals the
squared modulus of the normalized corner trace of \(g\).
For a free pair \(g,g'\), the operator \(T=1+g+g'\) is injective.
Indeed the matrix with diagonal entries zero and off-diagonal entries
\(T,T^*\) is the adjacency operator of the three-regular tree
(the universal cover of the graph with two vertices and three edges).
Its root resolvent is \(G(z)=(z-3m(z))^{-1}\), where the branch resolvent
\(m\) satisfies \(m(z)=(z-2m(z))^{-1}\) and \(m(z)\sim z^{-1}\).
These identities follow by taking the Schur complement at the root.
Since \(m(i\alpha)\to-i/\sqrt{2}\) as \(\alpha\downarrow0\), the
spectral measure has mass \(\lim_{\alpha\downarrow0}i\alpha G(i\alpha)=0\)
at zero. Thus the determinant on the kernel complement agrees with
the ordinary Fuglede--Kadison determinant. Ben Aribi's determinant formula
\cite[Theorem~0.2]{BenAribi} and multiplicativity of the
Fuglede--Kadison determinant \(\Delta\) yield
\[
 \Delta(B_r)=(2/\sqrt3)^r,\qquad \|B_r\|\le3^r .
\]
Fix \(1<\lambda<2/\sqrt3\). Integrate \(\log|B_r|\) on the two sides
of the cutoff \(\lambda^r\). The spectral projection for
\([\lambda^r,\infty)\) has normalized trace at least
\[
 q_0=\frac{\log(2/(\sqrt3\lambda))}{\log(3/\lambda)}>0
       \qquad (r\ge1).
\]
At \(s=q_0/4\) the product inequality for generalized singular values
\cite{FackKosaki} gives
\[
  \lambda^r\le\mu_{2s}(B_r)
  \le\mu_s(B_r\zeta_j)\,\mu_s(\zeta_j^{-1})
  \le(C_1/s)\mu_s(B_r\zeta_j).
\]
The coarse root corner has ambient trace \(a_j^2\).
By \eqref{q:terminal}, the terminal derivative has singular value
at normalized index \(s\) at least
\(t_j=(s/C_1)a_j^{-1/2}e^{-\sqrt h}\lambda^{H_j-1}\).
For large \(j\), \(\log(1+t_j)\ge c_1 H_j\) for a fixed \(c_1>0\).
Compression decreases generalized singular values. Thus, if
\(D=\dot{[X]}\), then for every \(\varepsilon>0\),
\[
 \tau_N\big((\log(1+|D|))^{1+\varepsilon}\big)
 \ \ge\ s a_j^2 (c_1H_j)^{1+\varepsilon}
 \ \ge\ c_\varepsilon d^{-4}2^{\,2\varepsilon d^2-20d}.
\]
Letting \(j\to\infty\) proves divergence at \((0,\Pi)\). Right
multiplication of every derivative by \(\Pi\) sends a solution at
\((0,1)\) to a solution at \((0,\Pi)\), and \(\mu_t(D'\Pi)\le\mu_t(D')\).
Finally, the row and column projections \(p_I,p_J\) for \([X]\) are
disjoint. Differentiating the support identities for \([X]\) and
\([X]^*\) gives \(L_{p_I}R_{p_J}\dot{(\operatorname{Re}[X])}=D/2\).
Consequently the single self-adjoint variable
\(\operatorname{Re}[X]\) has the asserted property.
\end{proof}

\section{A complex proof of the packing estimate at quadratic tolerance}
\label{cx:sec}

\begin{lemma}\label{cx:pack}
Let $\Phi\colon\R^a\times\R^m\to\R$ be a polynomial of degree at most
$p$, where $a\ge1$, $m\ge0$, and $d=a+m$. Write
$z=\nabla_x\Phi(x,u)$ and $v=\nabla_u\Phi(x,u)$. For $R\ge1$ and
$0<t<1$, put
\[
 Q_t=\{(x,z):\text{ for some }u,\quad
 |x|,|z|,|u|\le R,\quad |v|\le Rt^2\}.
\]
There is a constant $C=C(R,p)$ such that
$\Pack_t(Q_t)\le C^d t^{-a}$.
\end{lemma}

\begin{proof}
In the coordinates $(x,z,u,v)$ on $\C^{2d}$ use the complex symplectic
form $\Omega=\sum_i dx_i\wedge dz_i+\sum_j du_j\wedge dv_j$.
Suppose that $L\subset\C^{2d}$ is a complex Lagrangian plane, and
let $A$ be the matrix whose columns form an orthonormal basis of $L$.
For a set $I$ of $d$ coordinate indices, let $A_I$ be the corresponding
row minor and let $j(I)$ count the indices in the block $(x,z)$.
We claim that
\begin{equation}\label{cx:minors}
    \sum_{|I|=d,\ j(I)\le a}|\det A_I|^2\ge\tfrac12.
\end{equation}
On each symplectic pair let $K(q,h)=(-\overline h,\overline q)$,
and let $\sigma$ interchange the indices of each pair. Since
$\Omega|_L=0$, the antiunitary map $K$ takes $L$ onto $L^\perp$.
The matrix $(A,KA)$ is therefore unitary. Its complementary minors
give $|\det A_I|=|\det(KA)_{I^c}|=|\det A_{\sigma(I^c)}|$.
The involution $I\mapsto\sigma(I^c)$ takes $j(I)$ to $2a-j(I)$.
The Cauchy--Binet identity $\sum_{|I|=d}|\det A_I|^2=1$
now proves \eqref{cx:minors}.

Complexify $\Phi$ and apply the linear map
$(x,z,u,v)\mapsto(x,z,tu,v/t)$ to its complex gradient graph.
The resulting graph $G_t$ is a closed smooth complex submanifold of
complex dimension $d$. Its tangent planes are Lagrangian by symmetry
of the Hessian of $\Phi$ and invariance of $\Omega$ under the linear map.
Lift any $s$ points of a $t$-separated subset of $Q_t$ to real points
of $G_t$. The ambient balls centered at the lifts with radius $t/3$
are disjoint. Let $V$ denote the part of $G_t$ in their union and let
$\omega_k$ be the volume of the real unit $k$-ball, with $\omega_0=1$.
The complex manifold $G_t$ is calibrated by the $d$th divided power of
the Euclidean K\"ahler form \cite{HL}. The monotonicity formula for minimal
submanifolds \cite{Federer}, applied at each smooth center of density
one, gives
\begin{equation}\label{cx:lower}
     \operatorname{vol}_{2d}(V)\ge s\,\omega_{2d}(t/3)^{2d}.
\end{equation}

For the coordinate projection $\pi_I\colon G_t\to\C^d$, the
real Jacobian on a tangent plane with frame $A$ is $|\det A_I|^2$.
In the original position variables $(x,u)$, the map $\pi_I$ is
polynomial of degree at most $h=\max(1,p-1)$. Each regular fiber
has at most $h^d$ points by the affine B\'ezout bound for isolated
zeros \cite[Example~8.4.6]{Bezout}. The set of critical values has Lebesgue measure
zero. Integration of \eqref{cx:minors} and the area formula
\cite{Federer} therefore give
\[
   \operatorname{vol}_{2d}(V)
   \le 2h^d\!\sum_{|I|=d,\ j(I)\le a}
                        \operatorname{vol}_{2d}(\pi_I(V)).
\]
Put $B=3(R+1)$. At the centers of the balls the block $(x,z)$ is real
and the block $(tu,v/t)$ has norm at most $\sqrt{2}Rt$. Thus for
$j=j(I)$ the image $\pi_I(V)$ lies in the product of three real balls,
of dimensions $j,j,2(d-j)$ and radii $B,t/3,Bt$, respectively.
The product of the corresponding unit balls lies in the ball of radius
$\sqrt{3}$ in $\R^{2d}$. It follows that, whenever $j\le a$,
\[
   \operatorname{vol}_{2d}(\pi_I(V))
       \le 3^d\omega_{2d} B^{2d}t^{2d-j}
       \le 3^d\omega_{2d} B^{2d}t^{2d-a}.
\]
There are at most $4^d$ index sets $I$. Comparing this bound with
\eqref{cx:lower} gives $s\le (216h B^2)^d t^{-a}$.
\end{proof}

For the nonvanishing argument, the two auxiliary gradient blocks can
each be bounded by $\alpha=t^2/2$ before the matrix size is chosen.
Their joint norm is at most $\sqrt{2}\alpha\le t^2$, so
Lemma~\ref{cx:pack} gives the required upper packing bound.
The regularized determinant estimate uses the linear tolerance $Rt$
of Lemma~\ref{lem:linearpack}.


\begin{thebibliography}{99}
\raggedright

\bibitem{Alekseev}
V. Alekseev,
\href{https://doi.org/10.1016/j.jfa.2014.10.005}{\emph{On the first continuous \(L^2\)-cohomology of free group factors}},
J. Funct. Anal. \textbf{267} (2014), no.~11, 4264--4279.

\bibitem{Atiyah}
M.~F. Atiyah,
\href{https://numdam.org/item/AST_1976__32-33__43_0/}{\emph{Elliptic
operators, discrete groups and von Neumann algebras}},
Ast\'erisque \textbf{32--33} (1976), 43--72.

\bibitem{BenAribi}
F. Ben Aribi, \href{https://doi.org/10.5802/cml.79}{\emph{Fuglede--Kadison
determinants over free groups and Lehmer's constants}},
Confluentes Math. \textbf{14} (2022), no.~1, 3--22.

\bibitem{CS}
A. Connes and D. Shlyakhtenko,
\href{https://doi.org/10.1515/crll.2005.2005.586.125}{\emph{\(L^2\)-homology
for von Neumann algebras}}, J. Reine Angew. Math. \textbf{586} (2005), 125--168.

\bibitem{FackKosaki}
T. Fack and H. Kosaki,
\href{https://doi.org/10.2140/pjm.1986.123.269}{\emph{Generalized \(s\)-numbers
of \(\tau\)-measurable operators}}, Pacific J. Math. \textbf{123} (1986),
no.~2, 269--300.

\bibitem{Federer}
H. Federer, \href{https://doi.org/10.1007/978-3-642-62010-2}{\emph{Geometric
measure theory}}, Grundlehren Math. Wiss., vol.~153, Springer-Verlag,
New York, 1969.

\bibitem{Bezout}
W. Fulton, \href{https://doi.org/10.1007/978-1-4612-1700-8}{\emph{Intersection
theory}}, second ed., Ergeb. Math. Grenzgeb. (3), vol.~2, Springer-Verlag,
Berlin, 1998.

\bibitem{HT}
U. Haagerup and S. Thorbj{\o}rnsen,
\href{https://doi.org/10.4007/annals.2005.162.711}{\emph{A new application of
random matrices: \(\operatorname{Ext}(C^*_{\mathrm{red}}(F_2))\) is not a group}},
Ann. of Math. (2) \textbf{162} (2005), no.~2, 711--775.

\bibitem{HL}
R. Harvey and H.~B. Lawson, Jr.,
\href{https://doi.org/10.1007/BF02392726}{\emph{Calibrated geometries}},
Acta Math. \textbf{148} (1982), 47--157.

\bibitem{Jung}
K. Jung, \href{https://msp.org/pjm/2003/211-2/pjm-v211-n2-p04-p.pdf}{\emph{A
free entropy dimension lemma}}, Pacific J. Math. \textbf{211} (2003),
no.~2, 265--271.

\bibitem{Kechris}
A.~S. Kechris, \href{https://doi.org/10.1007/978-1-4612-4190-4}{\emph{Classical
descriptive set theory}}, Grad. Texts in Math., vol.~156, Springer-Verlag,
New York, 1995.

\bibitem{Luck}
W. L\"uck, \href{https://doi.org/10.1007/978-3-662-04687-6}{\emph{\(L^2\)-invariants:
theory and applications to geometry and \(K\)-theory}}, Ergeb. Math. Grenzgeb. (3),
vol.~44, Springer-Verlag, Berlin, 2002.

\bibitem{NS}
A. Nica and R. Speicher,
\href{https://doi.org/10.1353/ajm.1996.0034}{\emph{On the multiplication of
free \(N\)-tuples of noncommutative random variables}}, Amer. J. Math.
\textbf{118} (1996), no.~4, 799--837.

\bibitem{PV}
S. Popa and S. Vaes,
\href{https://doi.org/10.1093/imrn/rnu054}{\emph{Vanishing of the first
continuous \(L^2\)-cohomology for \(\mathrm{II}_1\) factors}},
Int. Math. Res. Not. IMRN (2015), no.~12, 3899--3907.

\bibitem{Shlyakhtenko2026}
D. Shlyakhtenko, \emph{On the isomorphism of semicircular algebras},
preprint, October 2026, \href{https://arxiv.org/abs/2610.19590}{arXiv:2610.19590}.

\bibitem{Thom}
A. Thom, \href{https://doi.org/10.1007/s00039-007-0634-7}{\emph{\(L^2\)-cohomology
for von Neumann algebras}}, Geom. Funct. Anal. \textbf{18} (2008),
no.~1, 251--270.

\bibitem{Voi}
D. Voiculescu, \href{https://doi.org/10.1007/BF01245072}{\emph{Limit laws for
random matrices and free products}}, Invent. Math. \textbf{104} (1991),
no.~1, 201--220.

\end{thebibliography}
\end{document}